\pdfoutput=1
\RequirePackage{ifpdf}
\ifpdf
\documentclass[pdftex]{sigma}
\else
\documentclass{sigma}
\fi

\newcommand{\id}{\, \mathrm{d}}
\newcommand{\ee}{\mathrm{e}}
\newcommand{\im}{\mathrm{i}}
\newcommand{\bO}{\mathcal{O}}

\begin{document}

\allowdisplaybreaks

\newcommand{\arXivNumber}{1609.02827}

\renewcommand{\PaperNumber}{101}

\FirstPageHeading

\ShortArticleName{Uniform Asymptotic Expansion for the Incomplete Beta Function}

\ArticleName{Uniform Asymptotic Expansion\\ for the Incomplete Beta Function}

\Author{Gerg\H{o} NEMES and Adri B.~OLDE DAALHUIS}

\AuthorNameForHeading{G.~Nemes and A.B.~Olde Daalhuis}
\Address{Maxwell Institute and School of Mathematics, The University of Edinburgh,\\
 Peter Guthrie Tait Road, Edinburgh EH9 3FD, UK}
\Email{\href{mailto:Gergo.Nemes@ed.ac.uk}{Gergo.Nemes@ed.ac.uk}, \href{mailto:A.B.Olde.Daalhuis@ed.ac.uk}{A.B.Olde.Daalhuis@ed.ac.uk}}
\URLaddress{\url{http://www.maths.ed.ac.uk/~gnemes/}, \url{http://www.maths.ed.ac.uk/~adri/}}

\ArticleDates{Received September 12, 2016, in f\/inal form October 21, 2016; Published online October 25, 2016}

\Abstract{In [Temme N.M., Special functions. An introduction to the classical functions of mathematical physics, \textit{A Wiley-Interscience Publication}, John Wiley \& Sons, Inc., New York, 1996, Section~11.3.3.1] a uniform asymptotic expansion for the incomplete beta function was derived. It was not obvious from those results that the expansion is actually an asymptotic expansion. We derive a remainder estimate that clearly shows that the result indeed has an asymptotic property, and we also give a recurrence relation for the coef\/f\/icients.}


\Classification{41A60; 33B20}

\section{Introduction}

For positive real numbers $a$, $b$ and $x\in[0,1]$, the (normalised) incomplete beta function $I_x(a,b)$ is def\/ined by
\begin{gather*}
I_x(a,b)=\frac{1}{B(a,b)}\int_0^x t^{a-1}(1-t)^{b-1}\id t,
\end{gather*}
where $B(a,b)$ denotes the ordinary beta function:
\begin{gather*}
B(a,b)=\int_0^1 t^{a-1}(1-t)^{b-1}\id t=\frac{\Gamma(a)\Gamma(b)}{\Gamma(a+b)}
\end{gather*}
(see, e.g., \cite[\href{http://dlmf.nist.gov/8.17.i}{Section~8.17(i)}]{NIST:DLMF}). In this paper, we will use the notation of \cite[\href{http://dlmf.nist.gov/8.18.ii}{Section~8.18(ii)}]{NIST:DLMF}.

The incomplete beta function plays an important role in statistics in connection with the beta distribution (see, for instance, \cite[pp.~210--275]{JKB95}). Large parameter asymptotic approximations are useful in these applications. For f\/ixed $x$ and $b$, one could use the asymptotic expansion
\begin{gather}\label{Ihyp}
I_x(a,b)=\frac{x^a(1-x)^{b-1}}{a B(a,b)}\,{}_2F_1\left(\begin{matrix} {1,1 - b} \\ {a + 1} \end{matrix};\frac{x}{x-1}\right)
\sim \frac{x^a(1-x)^{b-1}}{a B(a,b)} \sum_{n=0}^\infty \frac{\left(1-b\right)_n}{\left(a+1\right)_n}\left(\frac{x}{x-1}\right)^n,
\end{gather}
as $a\to+\infty$. The right-hand side of \eqref{Ihyp} converges only for $x\in [0,\frac12)$, but for any f\/ixed $x\in [0,1)$ it is still useful when used as an asymptotic expansion as $a\to+\infty$. For more details, see \cite[\href{http://dx.doi.org/10.1002/9781118032572.ch11}{Section~11.3.3}]{Temme96}. However, it is readily seen that~\eqref{Ihyp} breaks down as $x\to 1$. Since this limit has signif\/icant importance in applications, Temme derived in \cite[\href{http://dx.doi.org/10.1002/9781118032572.ch11}{Section~11.3.3.1}]{Temme96} an asymptotic expansion as $a\to+\infty$ that holds uniformly for $x\in (0,1]$. His result can be stated as follows.

\begin{theorem}
Let $\xi=-\ln x$. Then for any f\/ixed positive integer $N$ and fixed positive real $b$,
\begin{gather}\label{Thm}
I_x(a,b)=\frac{\Gamma(a+b)}{\Gamma(a)}\left(\sum_{n=0}^{N-1} d_n F_n+\bO\big(a^{-N}\big) F_0\right),
\end{gather}
as $a\to+\infty$, uniformly for $x\in(0,1]$. The functions $F_n=F_n(\xi,a,b)$ are defined by the recurrence relation
\begin{gather}\label{Frec}
aF_{n+1}=(n+b-a\xi)F_n+n\xi F_{n-1},
\end{gather}
with
\begin{gather*}
F_0=a^{-b}Q(b,a\xi),\qquad F_1=\frac{b-a\xi}{a}F_0+\frac{\xi^b \ee^{-a\xi}}{a\Gamma(b)},
\end{gather*}
and $Q(a,z)=\Gamma(a,z)/\Gamma(a)$ is the normalised incomplete gamma function $($see {\rm \cite[\href{http://dlmf.nist.gov/8.2.i}{Section~8.2(i)}]{NIST:DLMF})}. The coefficients $d_n=d_n(\xi,b)$ are defined by the generating function
\begin{gather}\label{Gdn}
\left(\frac{1-\ee^{-t}}{t}\right)^{b-1}=\sum_{n=0}^{\infty}d_{n}(t-\xi)^{n}.
\end{gather}
In particular,
\begin{gather*}
d_0=\left(\frac{1-x}{\xi}\right)^{b-1},\qquad d_1=\frac{x\xi +x-1}{(1-x)\xi}(b-1) d_0.
\end{gather*}
They satisfy the recurrence relation
\begin{gather}
\xi(n+1)(n+2)d_0d_{n+2} =\xi\sum_{m=0}^{n}(m+1)\left(n-2m+1+\frac{m-n-1}{b-1}\right)d_{m+1}d_{n-m+1} \nonumber\\
\hphantom{\xi(n+1)(n+2)d_0d_{n+2} =}{} +\sum_{m=0}^{n}(m+1)\left(n-2m-2-\xi+\frac{m-n}{b-1}\right)d_{m+1}d_{n-m} \nonumber\\
\hphantom{\xi(n+1)(n+2)d_0d_{n+2} =}{} +\sum_{m=0}^n (1-m-b)d_{m}d_{n-m}.\label{drec}
\end{gather}
In the case that $b=1$, we have $d_0=1$ and $d_n=0$ for $n \geq 1$.
\end{theorem}

Our contribution is the remainder estimate in \eqref{Thm} and the recurrence relation~\eqref{drec}. In fact, it is not at all obvious from~\eqref{Frec} that the sequence $\{F_n\}_{n=0}^{\infty}$ has an asymptotic property as $a\to+\infty$. We will show that for any non-negative integer~$n$,
\begin{gather}\label{Fasym}
0< F_{n+1}\leq \frac{n+\beta}{a}F_n,
\end{gather}
where $\beta=\max(1,b)$.

In \cite[\href{http://dx.doi.org/10.1142/9789814612166_0038}{Section~38.2.8}]{Temme15} the function $F_n$ is identif\/ied as a Kummer $U$-function:
\begin{gather*}
F_n=\frac{\xi^{n+b}\ee^{-a\xi}n!}{\Gamma(b)}U (n+1,n+b+1,a\xi ).
\end{gather*}

\section{Proof of the main results}\label{Proof}

We proceed similarly as in \cite[\href{http://dx.doi.org/10.1002/9781118032572.ch11}{Section~11.3.3.1}]{Temme96} and start with the integral representation
\begin{gather}\label{IxabInt}
I_x(a,b)=\frac{1}{B(a,b)}\int_\xi ^{+\infty} t^{b-1}\ee^{-at}\left(\frac{1-\ee^{-t}}{t}\right)^{b-1}\id t.
\end{gather}
We substitute the truncated Taylor series expansion
\begin{gather*}
\left(\frac{1-\ee^{-t}}{t}\right)^{b-1}=\sum_{n=0}^{N-1}d_{n}(t-\xi)^{n}+r_N(t)
\end{gather*}
into \eqref{IxabInt} and obtain
\begin{gather*}
I_x(a,b)=\frac{\Gamma(a+b)}{\Gamma(a)}\left(\sum_{n=0}^{N-1} d_n F_n+R_N(a,b,x)\right),
\end{gather*}
where $F_n$ is given by the integral representation
\begin{gather}\label{FInt}
F_n=\frac{1}{\Gamma(b)}\int_\xi ^{+\infty} t^{b-1}\ee^{-at}(t-\xi)^{n}\id t=\frac{\ee^{-a\xi}}{\Gamma(b)}\int_0 ^{+\infty} (\tau+\xi)^{b-1}\tau^n\ee^{-a\tau}\id \tau,
\end{gather}
and the remainder term $R_N(a,b,x)$ is def\/ined by
\begin{gather}\label{RN}
R_N(a,b,x)=\frac{1}{\Gamma(b)}\int_\xi ^{+\infty} t^{b-1}\ee^{-at}r_N(t)\id t.
\end{gather}
The recurrence relation \eqref{Frec} can be obtained from \eqref{FInt} via a simple integration by parts.

Let, for a moment,
\begin{gather*}
c_n(a,b)=\int_0 ^{+\infty} (\tau+\xi)^{b-1}\tau^n\ee^{-a\tau}\id \tau.
\end{gather*}
Then via integration by parts we f\/ind
\begin{gather}\label{cnRec}
ac_{n+1}(a,b)=(n+b)c_{n}(a,b)+\xi (1-b)c_{n}(a,b-1).
\end{gather}
We make the observation that
\begin{gather}\label{cnm}
0\leq\xi c_{n}(a,b-1)=\xi\int_0 ^{+\infty} (\tau+\xi)^{b-2}\tau^n\ee^{-a\tau}\id \tau\leq c_{n}(a,b).
\end{gather}
It follows from \eqref{cnRec} and \eqref{cnm} that
\begin{gather*}
ac_{n+1}(a,b)\leq
\begin{cases}
(n+1) c_n(a,b) & \text{if} \ 0<b\leq1, \\
(n+b) c_n(a,b) & \text{if} \ b\geq1.
\end{cases}
\end{gather*}
Since $F_n= \ee^{-a\xi} c_n(a,b)/\Gamma(b)$, this inequality implies~\eqref{Fasym}.

To obtain the remainder estimate in~\eqref{Thm}, we use the Cauchy integral representation
\begin{gather}\label{rNint}
r_N(t)=\frac{(t-\xi)^N}{2\pi\im}\oint_{\{\xi,t\}} \frac{\left(\frac{1-\ee^{-\tau}}{\tau}\right)^{b-1}}{(\tau-t)(\tau-\xi)^N}\id\tau,
\end{gather}
where the contour encircles the points $\xi$ and $t$ once in the positive sense. From the integral representation \eqref{RN}, we have that $0\leq\xi\leq t$. Thus, in the case that $N\geq1$, we can deform the contour in~\eqref{rNint} to the path
\begin{gather*}
[1+\infty\im,1+\pi\im]\cup [1+\pi\im,-1+\pi\im]\cup [-1+\pi\im,-1-\pi\im]\\
\hphantom{[1+\infty\im,1+\pi\im]}{} \cup [-1-\pi\im,1-\pi\im]\cup [1-\pi\im,1-\infty\im].
\end{gather*}
For the integrals along the f\/inal three portions of the path, we have the estimates
\begin{gather}
\left|\frac{1}{2\pi\im}\int_{-1+\pi\im}^{-1-\pi\im} \frac{\left(\frac{1-\ee^{-\tau}}{\tau}\right)^{b-1}}{(\tau-t)(\tau-\xi)^N}\id\tau\right|\leq
\frac{\max\left((\ee-1)^{b-1},\left(\frac{\ee+1}{\sqrt{\pi^2+1}}\right)^{b-1}\right)}{(1+\xi)^{N+1}},
\nonumber\\
\label{rNint2}
\left|\frac{1}{2\pi\im}\int_{-1-\pi\im}^{1-\pi\im} \frac{\left(\frac{1-\ee^{-\tau}}{\tau}\right)^{b-1}}{(\tau-t)(\tau-\xi)^N}\id\tau\right|\leq
\frac{\max\left(\left(\frac{\ee^{\pm1}+1}{\sqrt{\pi^2+1}}\right)^{b-1}\right)}{\pi^{N+2}},
\end{gather}
and
\begin{gather}
\left|\frac{1}{2\pi\im}\int_{1-\pi\im}^{1-\infty\im} \frac{\left(\frac{1-\ee^{-\tau}}{\tau}\right)^{b-1}}{(\tau-t)(\tau-\xi)^N}\id\tau\right| \leq
\frac{1}{2\pi}
\int_\pi^{+\infty}\frac{\max\big(\big(1\pm\ee^{-1}\big)^{b-1}\big) \big(s^2+1 \big)^{(1-b)/2}}{\sqrt{s^2+ (1-t )^2}\big(s^2+(1-\xi)^2\big)^{N/2}}\id s \nonumber \\
\hphantom{\left|\frac{1}{2\pi\im}\int_{1-\pi\im}^{1-\infty\im} \frac{\left(\frac{1-\ee^{-\tau}}{\tau}\right)^{b-1}}{(\tau-t)(\tau-\xi)^N}\id\tau\right| }{}
 \leq \frac{\max\big(\big(1\pm\ee^{-1}\big)^{b-1}\big)}{2\pi}
\int_\pi^{+\infty}\frac{ \big(s^2+1 \big)^{(1-b)/2}}{s^{N+1}}\id s,\label{rNint3}
\end{gather}
respectively. The integrals along the f\/irst two portions can be estimated similarly to~\eqref{rNint2} and~\eqref{rNint3}. Hence, for $0\leq\xi\leq t$ and $N\geq1$, we have
\begin{gather*}
|r_N(t)|\leq C_N(b) (t-\xi)^N,
\end{gather*}
where the constant $C_N(b)$ does not depend on $\xi$. Using this result in the integral representation~\eqref{RN}, we can infer that
\begin{gather*}
|R_N(a,b,x)|\leq C_N(b)F_N.
\end{gather*}
Finally, combining this result with the inequalities \eqref{Fasym}, we obtain the required remainder estimate in \eqref{Thm}.

The reader can check that the function $f(t)=\left(\frac{1-\ee^{-t}}{t}\right)^{b-1}$ is a solution of the nonlinear dif\/ferential equation
\begin{gather*}
tf(t)f''(t)-\frac{b-2}{b-1}tf'^2(t)+(t+2)f(t)f'(t)+(b-1)f^2(t)=0.
\end{gather*}
If we substitute the Taylor series \eqref{Gdn} into this dif\/ferential equation and rearrange the result, we obtain the recurrence relation~\eqref{drec}.

\subsection*{Acknowledgements}

This research was supported by a research grant (GRANT11863412/70NANB15H221) from the National Institute of Standards and Technology. The authors thank the anonymous referees for their helpful comments and suggestions on the manuscript.

\pdfbookmark[1]{References}{ref}
\LastPageEnding

\end{document}